\documentclass[12pt,twoside]{article}

\date{} 

\begin{document}

\centerline{} 

\centerline{} 

\centerline {\Large{\bf Tense operators on $m$--symmetric algebras}} 



\centerline{} 

\centerline{\bf {Aldo V. Figallo}} 

\centerline{} 

\centerline{Universidad Nacional de San Juan. Instituto de Ciencias B\'asicas} 

\centerline{Avda. I. de la Roza 230 (O), 5400 San Juan, Argentina} 

\centerline{and} 

\centerline{Universidad Nacional del Sur. Departamento de Matem\'atica} 

\centerline{Avda. Alem 1253, 8000 Bah\'ia Blanca, Argentina} 

\centerline{} 

\centerline{\bf {Carlos Gallardo}} 

\centerline{} 

\centerline{Universidad Nacional del Sur. Departamento de Matem\'atica} 

\centerline{Avda. Alem 1253, 8000 Bah\'ia Blanca, Argentina} 

\centerline{} 

\centerline{\bf {Gustavo Pelaitay}} 

\centerline{} 

\centerline{Universidad Nacional de San Juan. Instituto de Ciencias B\'asicas} 

\centerline{Avda. I. de la Roza 230 (O), 5400 San Juan, Argentina} 

\centerline{and} 

\centerline{Universidad Nacional del Sur. Departamento de Matem\'atica} 

\centerline{Avda. Alem 1253, 8000 Bah\'ia Blanca, Argentina} 

\centerline{gpelaitay@gmail.com}

\newtheorem{Theorem}{\quad Theorem}[section] 

\newtheorem{Definition}[Theorem]{\quad Definition} 

\newtheorem{Corollary}[Theorem]{\quad Corollary} 

\newtheorem{Lemma}[Theorem]{\quad Lemma} 

\newtheorem{Example}[Theorem]{\quad Example} 

\begin{abstract}Here we initiate an investigation of the equational classes of $m$--symmetric algebras endowed with two tense operators. These varieties is a generalization of tense algebras. Our main interest is the duality theory for these classes of algebras. In order to do this, we require Urquart's duality for Ockham algebras and Goldblatt's duality for bounded distributive lattice with operations. The dualities enable us to describe the lattices of congruences on tense $m$--symmetric algebras.
\end{abstract} 



\section{Introduction}

In 1977, generalizing De Morgan algebras by omitting the polarity condition (i.e.: the law of double negation), J. Berman (\cite{JB}) began the study of which he called distributive lattices with an additional unary operation.
Two years latter, A. Urquhart in \cite{U1} named them Ockham lattices with the justification that the so--called De Morgan laws are due, at least, in the case of propositional logic, to William  of Ockham. These algebras are the algebraic counterpart of logics provided with a negation operator which satisfies De Morgan laws.  Then recall that

\vspace{2mm}

An Ockham algebra is an algebra $\langle L, \wedge, \vee, N, 0,1\rangle$, where the reduct $\langle L, \wedge, \\ \vee, 0,1 \rangle$ is a bounded distributive lattice and $N$ is a unary operation satisfying the following conditions:\\ 

\begin{tabular}{ll}
(O1) $N(0) = 1$, &  \hspace{1.15cm} (O2) $N(1) = 0$,\\
(O3) $N(x \wedge y) = N(x) \vee N(y)$, & \hspace{1.3cm}(O4) $N(x \vee y)= N(x) \wedge N(y)$.
\end{tabular}

\vspace{3mm}

The name Ockham algebras has become classical and from  that moment on,  many articles have been published about this class of algebras. Many of the results obtained  have been reproduced in the important book by T. Blyth  and J. Varlet (\cite{TB.JV}), which  may be consulted by any reader interested in  broadening their knowledge on the topic.

\vspace{2mm}

For $m\geq 1$ and $n\geq 0$, $\mathcal{K}_{m,n}$ denotes the subvariety of Ockham algebras obtained by adjoining the equation $N^{2m+n}(x)=N^{n}(x)$. For $m\geq 1$ the algebras of the variety $\mathcal{K}_{m,0}$ are called $m$--symmetric algebras. (see \cite{TB.JV}, \cite{JVDZ}, \cite{GE})

The variety of $1$--symmetric algebras is the variety of De Morgan algebras which contains Kleene algebras and Boolean algebras as subvarieties. Kleene algebras are obtained by adjoining the equation $(x\wedge N(x))\vee (y\vee N(y))=y\vee N(y)$, whereas Boolean algebras are obtained by adjoining the equation $x\wedge N(x)=0$.

On the other hand, classical tense logic is a logical system obtained from bivalent logic by adding the tense operators $G$
({\it it is always going to be the case that}) and $H$ ({\it has always been the case that}). It is well--known that tense algebras represent the algebraic basis for bivalent tense logic \cite{B}, \cite{K}.

Starting with other logical systems and adding appropiate tense operators, we produce new tense logics (see \cite{DG,C}).

This paper deals with tense $m$--symmetric algebras, structures obtained from the $m$--symmetric algebras, by adding some tense operators.  
These algebras constitute a generalization of tense algebras. Our main interest is the duality theory
for these classes of algebras. In order to do this, we require Urquart's
duality for Ockham algebras and Goldblatt's duality for bounded distributive
lattice with operations. The dualities enable us to describe the
lattices of congruences on tense $m$--symmetric algebras.

\section{Preliminaries}\label{s2}

Even though the theory of Priestley spaces and its relation to bounded distributive lattices are well known (see \cite{P1}, \cite{P2} and \cite{P3}), we shall recall some definitions and results with the purpose of fixing the notations used in this paper.

Recall that a Priestley space (or $P$--space) is a compact totally disconnected ordered topological space. If $X$ is a $P$--space and $D(X)$ is the family of increasing, closed and open subsets of a $P$--space $X$, then $\langle D(X),\cap,\cup,\emptyset,X\rangle$ is a bounded distributive lattice.

On the other hand, let $L$ be a bounded distributive lattice and $X(L)$ be the set of all prime filters of $L$. Then $X(L)$ ordered by set inclusion and with the topology having as a sub--basis the sets $\sigma_{L}(a)=\{P\in X(L): a\in P\}$ and $X(L)\setminus \sigma_{L}(a)$ for each $a\in L$ is a $P$--space, and the mapping $\sigma_{L}:L\to D(X(L))$ is a lattice isomorphism. Besides, if $X$ is a $P$--space, the mapping $\varepsilon_{X}:X\to X(D(X))$ defined by $\varepsilon_{X}(x)=\{U\in D(X):x\in U\}$ is a homeomorphism and an order isomorphism.

If we denote by $\mathcal{L}$ the category of bounded distributive lattices and their corresponding homomorphisms and by $\mathcal{P}$ the category of $P$--spaces and the continuous increasing mappings (or $P$--functions), then there exists a duality between both categories by defining the contravariant functors $\Psi:\mathcal{P}\to\mathcal{L}$ and $\Phi:\mathcal{L}\to\mathcal{P}$ as follows:

\begin{itemize}
\item[(P1)] For each $P$--space $X$, $\Psi(X)=D(X)$ and for every $P$--function $f:X_1\to X_2$, $\Psi(f)(U)=f^{-1}(U)$ for all $U\in D(X_2)$.
\item[(P2)] For each bounded distributive lattice $L$, $\Phi(L)=X(L)$ and for every bounded lattice homomorphism $h:L_1\to L_2$, $\Phi(h)(F)=h^{-1}(F)$ for all $F\in X(L_2)$.
\end{itemize}

On the other hand, H. Priestley (\cite{P1,P2,P3}), proved that if $L$ is a bounded distributive lattice and $Y$ is a closed subset of $X(L)$, then

\begin{itemize}
\item[(P3)] $\Theta(Y)=\{(a,b)\in L\times L: \sigma_{L}(a)\cap Y=\sigma_{L}(b)\cap Y\}$ is a congruence on $L$ and that the correspondence $Y\mapsto \Theta(Y)$ is an anti--isomorphism from the lattice of all closed sets of $X(L)$ onto the lattice of all congruences on $L$.
\end{itemize}

The theory of Priestley topological duality was extended to Ockham algebras by A. Urquhart in \cite{U1}. This duality is equivalent to the one described now.

An Ockham space (or $OP$--space) is a pair $(X,g)$ where $X$ is a $P$--space and $g:X\to X$ is a decreasing continuous map. In addition, if $(X_1,g_1)$ and $(X_2,g_2)$ are $OP$--spaces, an $OP$--function from $(X_1,g_1)$ to $(X_2,g_2)$ is an increasing continuous function $f:X_1\to X_2$ such that $f\circ g_1=g_2\circ f$.

\vspace{2mm}

Then, the following results are fulfilled:

\begin{itemize}
\item[(U1)] If $(X,g)$ is an $OP$--space, then $(D(X),N_g)$ is an Ockham algebra where for all $U\in D(X)$, $N_g(U)=X\setminus g^{-1}(U)$,
\item[(U2)] If $(L,N)$ is an Ockham algebra, then $(X(L),g_N)$ is an $OP$--space where for all $P\in X(L)$, $g_{N}(P)=\{a\in L: N(a)\notin P\}$.  
\end{itemize}

Moreover, these constructions give a categorical dual equivalence.

\vspace{2mm}

On the other hand, in \cite{Gol} R. Goldblatt obtained a topological duality for bounded distributive lattices with operators, i.e. with a family of join--he\-mi\-mor\-phisms and/or meet--hemimorphisms. Now, we will describe this duality in the particular case of bounded distributive lattices endowed with two unary meet--hemimorphisms, $G$, $H$, which from now on will be called $G$--lattices.

\vspace{2mm}

A $gP$--space is a triple $(X,R_{G},R_{H})$ where $X$ is a $P$--space, $R_{G},R_{H}\subseteq X\times X$ are decreasing and the following conditions are satisfied:

\begin{itemize}
\item[{ (R1)}] for every $x\in X$, $R_{G}^{-1}(x)$ and $R_{H}^{-1}(x)$ are closed subsets of $X$,
\item[{ (R2)}] for each $U\in D(X)$, $G_{R_{G}}(U),H_{R_H}(U)\in D(X)$, where
\item[]$G_{R_{G}}(U)=\{y\in X: R_{G}^{-1}(y)\subseteq U\}$,  $H_{R_{H}}(U)=\{y\in X: R_{H}^{-1}(y)\subseteq U\}$.  
\end{itemize}

A $gP$--function from a $gP$--space $(X_1,R_{G_1},R_{H_1})$ into another one, $(X_2,R_{G_2},$ $R_{H_2})$, is a $P$--function $f:X_1\rightarrow X_2$ which satisfies the following conditions:

\begin{itemize}
\item[{ (r1)}] $(x,y)\in R_{G_1}$ implies $(f(x),f(y))\in R_{G_2}$ for $x,y\in X_1$,
\item[{ (r2)}] $(y,f(z))\in R_{G_2}$ implies that there is $x\in X_1$ such that $(x,z)\in R_{G_1}$ and $ f(x)\leq y$ for $z\in X_1$ and $y\in X_2$,
\item[{ (r3)}] $(x,y)\in R_{H_1}$ implies $(f(x),f(y))\in R_{H_2}$ for $x,y\in X_1$,
\item[{ (r4)}] $(y,f(z))\in R_{H_2}$ implies that there is $x\in X_1$ such that $(x,z)\in R_{H_1}$ and $f(x)\leq y$ for $z\in X_1$ and $y\in X_2$.   
\end{itemize}

In \cite{Gol} it was shown that
\begin{itemize}
\item[{ (G1)}] if $(L,G,H)$ is an $G$--lattice and $R^{L}_{T}\subseteq X(L)\times X(L)$ is defined by $R^{L}_{T}=\{(P,F)\in X(L)\times X(L): T^{-1}(F)\subseteq P \}$ for $T=G$ and $T=H$, then $(X(L),R_{G}^{L},R_{H}^{L})$ is a $gP$--space,  
\item[{ (G2)}] if $(X,R_{G},R_{H})$ is a $gP$--space, then $(D(X),G_{R_{G}},H_{R_{H}})$ is an $G$--lattice where $G_{R_{G}}$ and $H_{R_{H}}$ are defined in { (R2)}.
\end{itemize}

Taking into account { (G1)} and { (G2)} it is proved that the category of $gP$--spaces and $gP$--functions is dually equivalent to the category of $G$--lattices and their corresponding homomorphisms.

\section{Tense $m$--symmetric algebras}

In this section, we will describe a topological duality for tense $m$--symmetric algebras bearing in mind the results indicated in Section \ref{s2} and \cite{Gol}. 

\begin{Definition} A tense $m$--symmetric algebra is an algebra $\langle L,\vee,\wedge,N,G,$\\ $H,0,1\rangle$ such that the reduct $\langle L,\vee,\wedge,N,0,1\rangle$ is a bounded distributive lattice and $N$ is a dual endomorphism
satisfying the identity $N^{2m}(x)=x$ and 
$G$, $H$ are unary operators on $L$ verifying the following
conditions:
\begin{itemize}
  \item [{\rm (T1)}] $G(1)=1$, $H(1)=1$,
  \item [{\rm (T2)}] $G(x\wedge y)=G(x)\wedge G(y)$, $H(x\wedge y)=H(x)\wedge H(y)$,
  \item [{\rm (T3)}] $x\leq G(N( H(N^{2m-1} (x))))$, $x\leq H(N (G(N^{2m-1} (x))))$.
\end{itemize}
\end{Definition}

In what follows, we will denote these algebras by $(L,N,G,H)$ or simply by $L$ where no confusion may arise.

\vspace{2mm}

{\bf Remark} If $\langle L,\vee,\wedge,N,G,H,0,1\rangle$ is a tense $1$--symmetric algebra wich satisfies the identity $x\wedge Nx=0$, then $\langle L,\vee,\wedge,N,G,H,0,1\rangle$ is a tense algebra.

\vspace{2mm}

By $\textsf{TmS}$ we will denote the category of tense $m$--symmetric algebras and their corresponding homomorphisms.

\begin{Definition} A tense $m$--symmetric space (or $tms$--space) is a system $(X,g,R_{G},R_{H})$ where $(X,g)$ is an $OP$--space, $(X,R_{G},R_{H})$ is a $gP$--space and the following additional conditions are satisfied:
\begin{itemize}
  \item [{\rm (S1)}] $g^{2m}(x)=x$,
  \item [{\rm (S2)}] $(x,y)\in R_{G}$ implies $(g(y),g(x))\in R_{H}$,
  \item [{\rm (S3)}] $(x,y)\in R_{H}$ implies $(g(y),g(x))\in R_{G}$.
\end{itemize}
\end{Definition}

A $tms$--function from a $tms$--space $(X_1,g_1,R_{G_1},R_{H_1})$ into another one $(X_2,g_2,$\\ $R_{G_2},R_{H_2})$ is an $OP$--function $f:X_1\to X_2$ which is also a $gP$--function.

\vspace{2mm}

Let {\bf tmS} be the category whose objects are $tms$--space and whose morphisms are $tms$--functions.

\begin{Lemma}\label{p1} If $(X,g,R_{G},R_{H})$ is a $tms$--space, then $(D(X),N_{g},G_{R_G},H_{R_H})$ is a tense $m$--symmetric algebra, where for all $U\in D(X)$, $N_g(U)$, $G_{R_{G}}(U)$ and $H_{R_H}(U)$ are the sets defined in {\rm (U1)} and {\rm (R2)}, respectively.
\end{Lemma}

{\bf Proof.} We only prove {\rm (T3)}. Indeed, let $y\in U$ and $z\in R_{G}^{-1}(y)$. Suppose that $z\in g^{-1}(H_{R_H}(N_g^{2m-1}(U)))$. Then, $g(z)\in H_{R_H}(N_{g}^{2m-1}(U))$. From this last assertion, we infer that $R_{H}^{-1}(g(z))\subseteq N_{g}^{2m-1}(U)$. Besides, since $(z,y)\in R_{G}$, then by virtue of {\rm (S2)} we obtain that $g(y)\in R_{H}^{-1}(g(z))$ and so, $g(y)\in N_g^{2m-1}(U)$. Hence, $y\notin U$, which is a contradiction. Thus, $z\in N_g(H_{R_H}(N_g^{2m-1}(U)))$, from which we conclude that $R_{G}^{-1}(y)\subseteq N_{g}(H_{R_H}(N_{g}^{2m-1}($\\ $U)))$. So, $U\subseteq G_{R_{G}}(N_{g}(H_{R_H}(N_{g}^{2m-1}(U))))$. Similarly, it is proved that $U\subseteq H_{R_{H}}(N_{g}(G_{R_G}(N_{g}^{2m-1}(U))))$.

\vspace{1mm}

\begin{Lemma}\label{p2} If $(L,N,G,H)$ is a tense $m$--symmetric algebra, then $(X(L),$\\ $g_{N},R_{G}^{L},R_{H}^{L})$ is a $tms$--space, where $g_{N}$, $R_{G}^{L}$ and $R_{H}^{L}$ are those defined in {\rm (U2)} and {\rm (G1)}, respectively.
\end{Lemma}

{\bf Proof.} We only prove {\rm (S2)}. Indeed, let $a\in H^{-1}(g_{N}(P))$, and suppose that $a\notin g_{N}(F)$. Then $N(a)\in F$. Besides, from {\rm (T3)} we have that $N(a)\leq G(N(H(a)))$. So, $G(N(H(a)))\in F$. On the other hand, from hypothesis we infer that $G^{-1}(F)\subseteq P$. Hence, $N(H(a))\in P$. This last assertion allows us to conclude that $H(a)\notin g_{N}(P)$ which is a contradiction. Therefore, $a\in g_{N}(F)$ and so, $(g_{N}(F),g_{N}(P))\in R_{H}^{L}$.

\vspace{2mm}

From Lemma \ref{p1} and \ref{p2}, taking into account the results indicated in \cite{U1} and \cite{Gol}, we have Theorem \ref{t1}.

\begin{Theorem}\label{t1} The categories \textnormal{$\textsf{TmS}$} and {\bf tmS} are dually equivalent.
\end{Theorem}

This duality will allow us to determine the congruence lattice $Con_{tms}(L)$ of $L$.

\begin{Definition}\label{tsubset}
Let $(X,g,R_G,R_H)$ be a $tms$--space. A closed subset $Y$ of $X$ is a $tms$--subset of $X$ if it verifies these conditions for  $u,v\in X$:
\begin{itemize}
\item[\rm (tms1)] if $v\in R_T^{-1}(u)$ and $u\in Y$, then there is $w\in Y$ such that $w\in R_T^{-1}(u)$ and $w\leq v$, for $T=G$ and $T=H$.
\item[\rm (tms2)] $Y=g^{2m-1}(Y)$.
\end{itemize}
\end{Definition}

We will denote by $\mathcal{C}_{tms}(X)$ the set of all $tms$--subset of $X$.

\newpage

\begin{Lemma}\label{l3} Let $(L,N,G,H)$ be a tense $m$--symmetric algebra and $Y\in\mathcal{C}_{tms}(X(L))$. Then $\Theta(Y)$ is an $tms$--congruence, where $\Theta(Y)$ is defined as in {\rm (P3)}.
\end{Lemma}

{\bf Proof.} We only prove that $\Theta(Y)$ preserves $N$, $G$ and $H$. Let $(a,b)\in \Theta(Y)$. Since $Y=g^{2m-1}(Y)$, we have that $(Na,Nb)\in \Theta(Y)$. Indeed, suppose that $Q\in \sigma_{L}(N(a))\cap Y$. Then, $N(a)\in Q$ and $Q\in g^{2m-1}_{N}(Y)$. This last assertion allows us to infer that there is $P\in Y$ such that $g_{N}^{2m-1}(P)=Q$. So, since $N(a)\in g_{N}^{2m-1}(P)$ we have that $a\notin P$. Hence, $b\notin P$ and therefore $N(b)\in Q$. The other inclusion is proved in a similar way.

On the other hand, let $F\in G_{R_G}(\sigma_L(a))\cap Y$. Hence, $R_{G}^{-1}(F)\subseteq \sigma_{L}(a)$ and $F\in Y$. Suppose that $Q\in R_{G}^{-1}(F)$. Then, from (tms1) there is $W\in Y$ such that $W\subseteq Q$ and  $W\in R^{-1}_{G}(F)$. This last assertion allows us to infer that $W\in\sigma_{L}(a)$, from which we conclude that $W\in \sigma_{L}(b)\cap Y$. So, since $W\subseteq Q$ we have that $Q\in\sigma_{L}(b)$. Hence, $F\in G_{R_G}(\sigma_{L}(b))\cap Y$ and therefore, $G_{R_{G}}(\sigma_{L}(a))\cap Y\subseteq G_{R_{G}}(\sigma_{L}(b))\cap Y$. The other inclusion is proved in a similar way. Analogously, $\Theta(Y)$ preserves $H$.

\begin{Lemma}\label{l4} Let $(L,N,G,H)$ be a tense $m$--symmetric algebra, $\theta\in Con_{tms}(L)$ and $Y=\{\Phi(q)(F): F\in X(L/\theta)\}$, where $\Phi$ is defined as in {\rm (P2)}. Then, $Y\in\mathcal{C}_{tms}(X(L))$.
\end{Lemma}

{\bf Proof.} Let $\Phi(q):X(L/\theta)\to X(L)$ be the function defined by $\Phi(q)(F)=q^{-1}(F)$ for $F\in X(L/\theta)$. Then, since $Con_{tms}(L)$ is a sublattice of $Con(L)$ we have that $Y=\{\Phi(q)(F):F\in X(L/\theta)\}$ is a closed subset of $X(L)$ and $\theta=\Theta (Y)$. Besides, from \cite{Gol},\cite{U1} we have that $\Phi(q)$ is a $tms$--function. In addition, $Y$ is a $tms$--subset of $X(L)$. Indeed, let $V\in {(R^{L}_{G})}^{-1}(U)$ with $U\in Y$. From this last assertion, there is $F\in X(L/\theta)$ such that $\Phi(q)(F)=U$. Then, $(V,\Phi(q)(F))\in R_{G}^{L}$. So, from { (r2)} we have that there is $P\in X(L/\theta)$ such that $(P,F)\in R_{G}^{L/\theta}$ and $\Phi(q)(P)\subseteq V$. Therefore, from { (r1)} we infer that $\Phi(q)(P)\subseteq {(R_{G}^{L})}^{-1}(U)$. Similarly, it is proved (tms1) for $T=H$.

On the other hand, let $Q\in Y$. Then, there is $F\in X(L/\theta)$ such that $\Phi(q)(F)=Q$. Since, $\Phi(q)$ is an $OP$--function we have that $g_{N}(Q)=\Phi(q)(g_{N}(F))$. This last assertion and (S1) allows us to conclude that $Q\in g_{N}^{2m-1}(Y)$. Therefore, $Y\subseteq g_{N}^{2m-1}(Y)$. The other inclusion is straightforward. This completes the proof.

\vspace{2mm}

From Lemma \ref{l3} and \ref{l4}, we have Theorem \ref{t2}.

\begin{Theorem}\label{t2} Let $(L,N,G,H)$ be a tense $m$--symmetric algebra and $(X(L),$\\ $g_{N},R^{L}_{G},R^{L}_{H})$ be the associated $tms$--space of $L$. Then, there is an anti--iso\-morphisms between $Con_{tms}(L)$ and the lattice $\mathcal{C}_{tms}(X(L))$.
\end{Theorem}



\begin{thebibliography}{99} 



\bibitem{JB} {J. Bermann}, {\em Distributive Lattices with an additional unary operation}, Aequationes Math. {16} (1977), 165--171.

\bibitem{TB.JV} {T. Blyth and J. Varlet}, {\em Ockham Algebras}, Oxford University Press, New York, 1994. 

\bibitem{B} {J. Burges}, {\em Basic tense logic}, Handbook of philosophical logic, Vol. II, 89--133, Synthese Lib., 165, Reidel, Dordrecht, 1984.

\bibitem{C} {C. Chirita}, {\em Tense $\theta$--valued Moisil propositional logic}, Int. J. of Computers, Communications \& Control. {5} (2010),5, 642--653.

\bibitem{DG} {D. Diaconescu and G. Georgescu}, {\em Tense operators on $MV$--algebras and \L ukasiewicz--Moisil algebras}, Fund. Inform. {81} (2007),4, 379--408.

\bibitem{GE} {P. Garc\'ia and F. Esteva}, {\em On Ockham algebras: congruence lattices and subdirectly irreducible algebras}, Studia Logica, {55} (1995), 2, 319--346.

\bibitem{K} {T. Kowalski}, {\em Varieties of tense algebras}, Rep. Math. Logic. {32} (1998), 53--95. 

\bibitem{Gol} {R. Goldblatt}, {\em Varieties of complex algebras}, Ann. Pure Appl. Logic. {44} (1989),3, 173--242. 

\bibitem{P1} {H. Priestley}, {\em Representation of distributive lattices by means of ordered Stone spaces},
            Bull. London Math. Soc. { 2}(1970), 186--190. 

\bibitem{P2} {H. Priestley}, {\em Ordered topological spaces and the representation of distributive lattices},
            Proc. London Math. Soc. {24} (1972),3, 507--530. 

\bibitem{P3} {H. Priestley}, {\em Ordered sets and duality for distributive lattices}, Orders: description and roles (L'Arbresle, 1982), 39--60, North-Holland Math. Stud., 99, North-Holland, Amsterdam, 1984
\bibitem{U1} {A. Urquhart}, {\em Distributive lattices with a dual homomorphic operation}, Studia Logica {38} (1979), 201--209.

\bibitem{JVDZ} {J. Vaz De Carvalho}, {\em Congruences on algebras of ${\mathcal K}_{m,0}$}, Bull. Soc. Roy. Sci. Li\`ege { 54} (1985), 301--303.

\end{thebibliography}
\end{document}